\pdfoutput=1
\documentclass[12pt,reqno]{amsart}

\usepackage[utf8]{inputenc}
\usepackage[T1]{fontenc}
\usepackage[charter]{mathdesign}
\usepackage[draft=false,kerning=true]{microtype}

\usepackage{fullpage}

\usepackage{amsthm, amsmath, amsfonts}
\usepackage{mathtools}
\usepackage{enumerate}

\usepackage{graphicx}
\usepackage{tikz}
\usetikzlibrary{arrows,backgrounds,calc}
\usepackage{pgfplots}
\pgfplotsset{compat=1.12}

\usepackage{hyperref}

\newtheorem{theorem}{Theorem}

\theoremstyle{remark}
\newtheorem*{remark}{Remark}

\theoremstyle{definition}

\newtheorem*{example*}{Example}

\renewcommand{\MR}[1]{}

\DeclareSymbolFont{cmcal}{OMS}{cmsy}{m}{n}
\SetSymbolFont{cmcal}{bold}{OMS}{cmsy}{b}{n}
\DeclareSymbolFontAlphabet{\mathcal}{cmcal}

\DeclarePairedDelimiter{\abs}{\lvert}{\rvert}

\title[]{A Combinatorial Identity for Rooted Labeled Forests}

\author[B.~Hackl]{Benjamin Hackl}

\address[Benjamin Hackl]
{Institut f\"ur Mathematik,
  Alpen-Adria-Uni\-ver\-si\-t\"at Klagenfurt, Universit\"atsstra\ss e
  65--67, 9020 Klagenfurt, Austria}
\email{\href{mailto:benjamin.hackl@aau.at}{benjamin.hackl@aau.at}}
\thanks{B.~Hackl is supported by the Austrian
  Science Fund (FWF): P~28466-N35 and by a
  distinction grant from the Austrian Federal Ministry 
  of Education, Science and Research.}

\keywords{Combinatorial Identity, Forest, Set partition,
Hurwitz multinomial identity}
\subjclass[2010]{05A19; 05C05}

\begin{document}
\begin{abstract}
  In this brief note a straightforward combinatorial proof for
  an identity directly connecting rooted forests and unordered
  set partitions is provided. Furthermore, references that
  put this type of identity in the context of \emph{forest volumes}
  and \emph{multinomial identities} are given.
\end{abstract}

\maketitle

\section{The Identity}

The aim of this note is to provide an elementary and purely combinatorial proof
for an identity stated and proved (via induction) by 
Dorlas, Rebenko and Savoie in~\cite{Dorlas-Rebenko-Savoie:2019:forest-identity}.

Let $m$ and $p$ be positive integers with $p\leq m$ and define $\Omega = \{x_1, x_2, \ldots, x_m\}$ to be a set of $m$ variables. The
combinatorial identity of interest is
\begin{equation}\label{eq:identity}
  \sum_{P \in \Pi_p(\Omega)} \prod_{T\in P} 
  \Big(\sum_{x\in T} x \Big)^{\abs{T}-1}
   = \binom{m-1}{p-1} \Big(\sum_{x\in\Omega} x\Big)^{m-p},
\end{equation}
where $\Pi_p(\Omega)$ is the set of all set partitions of $\Omega$
that consist of $p$ parts.
The identity is particularly remarkable in the sense that it provides
a not immediately obvious connection between the ``partition world''
on the left-hand side and the ``binomial world'' on the right-hand 
side.

It should be noted that variations of~\eqref{eq:identity} are
well-known and can be found in, e.g., \cite[Theorem
5.3.4]{Stanley:1999:enum-comb-2}. Equations of this type are related
to Hurwitz multinomial identities, see \cite{Hurwitz:1902:abel-binomial},
\cite{Pitman:2002:forest-volume-decomposition} and 
\cite{Robertson:1964:abel-sum}. In particular, in 
\cite{Pitman:2002:forest-volume-decomposition} Pitman presents a 
systematic 
approach to interpreting identities of this type as decompositions of 
forest volumes, i.e., polynomials that enumerate special classes of
rooted forests. It is not too difficult to apply this framework
in order to prove~\eqref{eq:identity}---however, we want to present
a more explicit combinatorial proof based on double counting.

It is also easy to see that looking for a combinatorial interpretation
in the context of words (i.e., in a setting where the variables
in $\Omega$ are non-commutative as the multiplication is replaced
by concatenation) is not possible in general, as the following example
illustrates.

\begin{example*}
Take $m = 5$, $p = 2$ and consider the word $x_1 x_2 x_1$. On the 
right-hand side the word is constructed once within the parenthesis and
then enumerated by the binomial coefficient $\binom{m-1}{p-1} = 4$. 
However, on the left-hand
side the word can only be constructed from a single partition part
of size 4 containing both $x_1$ and $x_2$, i.e., from the partitions
\[
  \{x_1, x_2, x_3, x_4\}, \{x_5\}, \qquad 
  \{x_1, x_2, x_3, x_5\}, \{x_4\}, \qquad 
  \{x_1, x_2, x_4, x_5\}, \{x_3\}.
\]
Thus, the word $x_{1} x_{2} x_{1}$ occurs only three times on the
left-hand side. In
general, the left-hand side actually enumerates sets of words over
disjoint alphabets, as we take the product over an unordered set 
partition---in this particular setting, however, the word 
$x_1 x_2 x_1$ has to be constructed entirely from one of the parts of the
partition.
\end{example*}

\section{A Refined Version and its Interpretation}

We choose to prove a variation of~\eqref{eq:identity} that belongs to
the same family of identities as discussed
in~\cite[Section~4 and Corollary~8 in particular]{Pitman:2002:forest-volume-decomposition}.
\begin{theorem}
  Let $u$ be an additional variable. Then the identity
  \begin{equation}\label{eq:identity:u}
    \sum_{P \in \Pi(\Omega)} \prod_{T\in P} 
    u\Big(\sum_{x\in T} x \Big)^{\abs{T}-1}
    = u \Big(u +\sum_{x\in\Omega} x\Big)^{m-1}
  \end{equation}
  holds, where $\Pi(\Omega)$ is the set of set partitions of $\Omega$.
\end{theorem}
\begin{remark}
  Extracting the coefficient of $u^p$ on both sides 
  of~\eqref{eq:identity:u} immediately yields~\eqref{eq:identity}.
  As we will see in the proof, combinatorially, this identity
  describes two ways of constructing rooted forests on the
  vertex set $\{1,2,\ldots,m\}$: As a set of trees on a partition
  of the vertex set on the left-hand side, and as the result of
  deleting the root node in a rooted tree on an extended vertex set
  on the right-hand side.
\end{remark}
\begin{proof}
  We claim that both sides of~\eqref{eq:identity:u} enumerate rooted
  labeled forests on the vertex set $\{1,2,\ldots, m\}$ with respect
  to the number of components (enumerated
  by the variable $u$) and with respect to the out-degree of the
  vertices (enumerated by the variables $x_j$, respectively).
  
  We begin our proof by observing that as a consequence of the
  well-known Prüfer bijection (see 
  \cite[Theorem~5.19]{Bona:2007:introduction-enum}), the 
  multivariate ordinary generating function
  \begin{align}\notag
    T(u, x_1, x_2, \ldots, x_m) & = 
    (u + x_1 + x_2 + \dots + x_m)^{m-1} u x_1 x_2 \dots x_m \\ \label{eq:T:rhs}
    &=
    \Big(u + \sum_{x\in\Omega} x\Big)^{m-1} u \prod_{x\in\Omega} x
  \end{align}
  enumerates trees with vertex set $\{0, 1, \ldots, m\}$ where every 
  variable is associated to some vertex ($u$ is associated to $0$, 
  $x_j$ is associated to $j$ for all $1\leq j\leq m$) and 
  keeps track of the number of vertices adjacent to ``their'' vertex.
  In other words, this means that the
  coefficient of the monomial $u^p x_1^{d_1} x_2^{d_2} \dots x_m^{d_m}$
  in $T(u, x_1, x_2, \ldots, x_m)$ is precisely the number of trees 
  with $\deg(0) = p$ and $\deg(j) = d_j$ for all $1\leq j\leq m$.
  More details and a proof of this fact are given in, e.g., 
  \cite[Theorem~5.19]{Bona:2007:introduction-enum}.
  
  We are allowed to think of the trees enumerated by 
  $T(u, x_1, \ldots, x_m)$ as trees rooted at the vertex $0$. By
  removing the product $\prod_{x\in\Omega} x$ from~\eqref{eq:T:rhs}
  we obtain the right-hand side of our identity~\eqref{eq:identity:u}.
  Combinatorially, this corresponds to ignoring
  one neighbor for every vertex except for $0$. Equivalently, this can 
  be seen as keeping track of the out-degrees of the vertices instead of
  their degree. This is because when considering
  trees rooted at $0$, every vertex except $0$ has a parent in the
  tree.
  
  Finally, observe that by deleting vertex $0$, trees on the vertex 
  set $\{0,1,\ldots,m\}$ that are rooted at $0$ are in bijection
  to rooted forests on $\{1,2,\ldots,m\}$ where the out-degrees of
  the vertices $\{1,2,\ldots,m\}$ remain unchanged and the number of
  components in the forest corresponds to the degree of $0$. See
  Figure~\ref{fig:forest-tree-bijection} for an illustration.
  This proves our claim for the right-hand side
  of~\eqref{eq:identity:u}.
  
  \begin{figure}[ht]
  \centering
  \begin{tikzpicture}[edge from parent/.style={draw, -latex}, inner sep=1pt, font=\footnotesize]
    \draw node[draw, circle]{0}
    child{
      node[draw, circle]{1} edge from parent[dashed]
      child[solid]{
        node[draw, circle]{4}
        child[solid]{
          node[draw, circle]{3}
        }
        child[solid]{
          node[draw, circle]{6}
        }
      }
    }
    child{
      node[draw, circle]{5} edge from parent[dashed]
    }
    child{
      node[draw, circle]{8} edge from parent[dashed]
      child[solid] {node[draw, circle]{2}}
      child[solid] {node[draw, circle]{7}}
    };
  \end{tikzpicture}
  \caption{A tree with vertex set $\{0, 1, \ldots, 8\}$ rooted at $0$
  is in bijection to a rooted forest on $\{1,2,\ldots, 8\}$.}
  \label{fig:forest-tree-bijection}
  \end{figure}
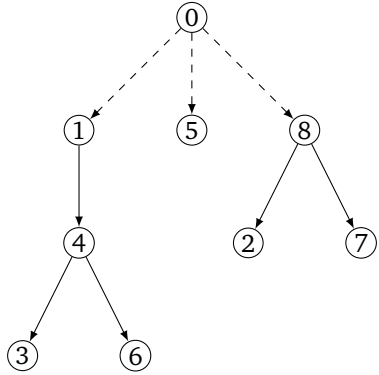

  The combinatorial interpretation of the left-hand side is
  quite straightforward. Any rooted forest on 
  $\{1,2,\ldots,m\}$ can be constructed by starting with an
  arbitrary set partition $P\in \Pi(\Omega)$, constructing
  rooted trees whose vertex sets are the subsets of $\{1,2,\ldots, m\}$
  corresponding to the partition parts $T\in P$, and combining
  these trees to a forest.
  
  Obviously, for every set partition $P$ the variable $u$ keeps 
  track of the number of parts in $P$---and, equivalently, of the 
  number of components of the forest. Finally, note that 
  for $\abs{T} \geq 2$ we can rewrite
  \[ \Big(\sum_{x\in T} x\Big)^{\abs{T}-1} 
  = \Big(\sum_{x\in T} x\Big)\Big(\sum_{x\in T} x\Big)^{\abs{T}-2}
  = \sum_{r\in T} r\Big(\sum_{x\in T} x\Big)^{\abs{T}-2}, \]
  which, by the same reasoning as above, enumerates all rooted trees
  on the vertices corresponding to $T$ with respect to the 
  out-degrees of the vertices. The same is true for $\abs{T} = 1$.
  
  Overall, both sides of~\eqref{eq:identity:u} can be interpreted as
  ordinary generating functions for the same family of combinatorial
  objects---which proves that the identity holds.
\end{proof}

\subsection*{Acknowledgements}
This note has been written during a research stay at Stellenbosch
University. I explicitly want to thank Helmut Prodinger for pointing
me towards the appropriate literature.

\bibliographystyle{amsplainurl}
\bibliography{forest-identity}

\end{document}